\documentclass[review]{elsarticle}
\usepackage{bm, booktabs, enumerate, graphicx, multirow, multicol, amsmath, amstext, amssymb, pdflscape, multirow, amsfonts, color, placeins, algorithm, algcompatible, setspace, hyperref, tikz, pgfplots, amsthm, rotating, caption, threeparttable, pgfgantt}
\usepackage{wrapfig,comment,tikz,animate,algpseudocode}

\usepackage{mathtools}
\DeclarePairedDelimiter{\ceil}{\lceil}{\rceil}

\usepackage[margin=0.9in]{geometry}
\usetikzlibrary{arrows}
\usepgfplotslibrary{groupplots}
\pgfplotsset{compat=1.16}

\biboptions{authoryear}
\allowdisplaybreaks

\begin{document}

\begin{frontmatter}
\title{Scheduling chemotherapy appointments under uncertainty by considering different nursing care delivery schemes}

\author{Serhat Gul\footnote{Corresponding author, e-mail: serhat.gul@tedu.edu.tr}}
\address{Department of Industrial Engineering, TED University, Ankara, Turkey}

\begin{abstract}
The flexibility level allowed in nursing care delivery and uncertainty in infusion durations are very important factors to be considered during the chemotherapy schedule generation task. The nursing care delivery scheme employed in an outpatient chemotherapy clinic (OCC) determines the strictness of the patient-to-nurse assignment policies, while the estimation of infusion durations affects the trade-off between patient waiting time and nurse overtime. We study the problem of daily scheduling of patients, assignment of patients to nurses and chairs under uncertainty in infusion durations for an OCC that functions according to any of the three commonly used nursing care delivery models representing fully flexible, partially flexible, and inflexible care models, respectively. We develop a two-stage stochastic mixed-integer programming model that is valid for the three care delivery models to minimize expected weighted cost of patient waiting time and nurse overtime. We propose multiple variants of a scenario grouping-based decomposition algorithm to solve the model using data of a major university oncology hospital. The variants of the algorithm differ from each other according to the method used to group scenarios. We compare input-based, solution-based and random scenario grouping methods within the decomposition algorithm. We obtain near-optimal schedules that are also significantly better than the schedules generated based on the policy used in the clinic. We analyze the impact of nursing care flexibility to determine whether partial or fully flexible delivery system is necessary to adequately improve waiting time and overtime. We examine the sensitivity of the performance measures to the cost coefficients and the number of nurses and chairs. Finally, we provide an estimation of the value of stochastic solution.       

\end{abstract}
\end{frontmatter}

\section{Introduction}	
The incidence rate of cancer and demand for chemotherapy services are expected to significantly rise by 2040 \citep{wilsonetal19}. This implies that the need for the efficient use of outpatient chemotherapy center (OCC) resources (i.e. nurses and chairs) will be even more critical than it is now. Therefore, the OCC managers has to face the challenge of developing intelligent approaches to create appointment schedules satisfying both providers and patients.     

In general, chemotherapy schedules are generated in two steps \citep{turkcanetal12, karakayaetal22}. The first step, called \textit{chemotherapy planning}, includes the allocation of patients to days, and followed by the second step, \textit{chemotherapy scheduling}. Patient appointment times are set and the nurses responsible for each chemotherapy treatment are determined in the second step. The current study focuses on the chemotherapy scheduling step.     

The two activities in the OCC that must be carefully modelled while designing schedules are premedication and infusion. Both activities are conducted when the patient is seated on a chemotherapy chair. In the premedication phase, the nurse assigned to a patient prepares them for infusion and gives premedication drugs to protect the patient from the side effects of chemotherapy drugs. A nurse can provide premedication service to only one patient at a time. However, they can handle multiple infusions simultaneously. At the infusion phase, the chemotherapy drugs flow from an intravenous bag to a catheter located in the veins of the patient according to the dosage determined by an oncologist. 

Coordinating the activities in the OCC through schedules while aiming to use nurses and chairs efficiently is a challenging task under uncertainty in infusion durations. When the system is totally flexible in terms of nursing care delivery scheme, relatively efficient schedules can be generated by taking into account the trade-off between patient and provider related metrics. However, when the patient-to-nurse assignment policies are made very strict with the intention of maximizing care quality, the resulting schedules may lead to low level of efficiency measures including patient waiting time and nurse overtime. Therefore, the level of flexibility for nursing care allowed in the OCC plays a critical role during the process of schedule generation. 

The fully-flexible system at which patients can be assigned to any one of the nurses is defined in the literature as \textit{functional care delivery model} \citep{liangturkcan16}. On the other hand, the \textit{primary care delivery model} represents the setting where a \textit{primary nurse} of a patient is responsible for the treatment of the patient at each chemotherapy visit to maintain care continuity. These two extreme settings are the most popular among clinics. In particular, 40\% of the clinics use functional care delivery model, while 39\% prefer primary care delivery models \citep{irelandetal04, liangturkcan16}. Alternative settings that blend the characteristics of these two settings are also implemented to improve care delivery by benefiting from the advantages of them. For example, in a care model which we call as \textit{partially flexible care delivery}, each patient may be assigned to their \textit{primary nurse} or to an \textit{alternative nurse} to some extent. Under this setting, the OCC manager gives priority to the assignment of patients to their primary nurses, but assign some patients to alternative nurses when necessary. In some settings, alternative nurses are further categorized as \textit{secondary nurses} and \textit{floating nurses} \citep{tabrizietal20}. In such a setting, each patient is also associated with a secondary nurse who is given priority over floating nurses during the nurse assignment. 

Nurse overtime and patient waiting time are the commonly considered conflicting metrics by the OCC managers to assess the schedules. To reduce operating costs of the OCC and improve job satisfaction levels of nurses, overtime must be reduced. The satisfaction levels of patients can be increased by reducing patient waiting times. Decreasing both waiting time and overtime by creating better schedules is a challenging task due to uncertainty in infusion durations \citep{karakayaetal22}. Estimating infusion duration for a patient may be difficult because of possible side effects that affect patients when chemotherapy drugs damage healthy cells. Terminating the treatment early, slowing down the administration, or changing medication plan may lead to deviations from the anticipated infusion durations. When the durations are underestimated while creating the daily chemotherapy schedules, patient waiting time may increase. On the other hand, overestimated durations may increase nurse overtime.

In this study, we consider the problem of sequencing patients, setting appointment times, and assigning patients to nurses and chairs under uncertainty in infusion durations for an OCC which may operate according to any of the commonly used care delivery models such as primary care, functional care or partially flexible care. We develop a two-stage stochastic mixed-integer programming (TSMIP) model to create a daily schedule by assuming that the list of patients to be treated on a day is already determined. We minimize the expected weighted cost of patient waiting time and nurse overtime by sampling scenarios from appropriate infusion duration distributions. The duration distributions are obtained from data of an OCC at a major university oncology hospital. We develop and test multiple versions of a scenario grouping-based decomposition (SGBD) algorithm on the TSMIP model. The algorithms differ from each other based on the scenario-grouping approach. We compare input-based, solution-based and random grouping approaches while implementing the SGBD algorithm. We compare the selected variant of the SGBD algorithm with a commercial solver to assess the optimality gap. We show the level of improvement provided by the near-optimal schedules with reference to the schedules obtained according to the policy used in the OCC. We analyze the impact of nursing care flexibility to determine whether a partial or fully flexible delivery system would be necessary to sufficiently improve the performance measures of the study. We conduct sensitivity analysis by examining the impact of the cost coefficients of waiting time and overtime, and the number of nurses and chairs to the performance measures. Finally, we provide an estimation of the value of stochastic solution.       

The TSMIP model in the current study is developed by extending the model in \citet{demiretal21}. We summarize the differences between the two studies next. \citet{demiretal21} provide a formulation which is applicable to only functional care delivery scheme. Since the focus is the nursing care flexibility in the current study, our model is applicable to primary care, partially flexible care delivery as well as the functional care delivery setting. To the best of our knowledge, the current study is the first that formulate a stochastic appointment scheduling model which is valid for all these three chemotherapy delivery settings. Note that the nurse and chair assignment decisions are considered in the second stage of the stochastic programming (SP) model in \citet{demiretal21}. On the other hand, they are considered in the first stage of our model, which makes the SP model structure very different. We do not consider idle time in the model, while \citet{demiretal21} does, because they also show that reducing overtime leads to reduced idle time. Finally, they solve the model using a progressive hedging algorithm. We solve the TSMIP model by implementing an SGBD algorithm with its multiple variants. 

In the next section, we provide a detailed literature review on chemotherapy scheduling. Next, in Section 3, we define the problem and formulate the TSMIP model. In Section 4, we give the details of the SGBD algorithms we propose. In Section 5, we conduct a comprehensive numerical experiments and discuss their results. In Section 6, we provide conclusions and briefly discuss potential future studies. 

\section{Literature Review}
The literature review consists of two categories of articles. The first category includes deterministic chemotherapy scheduling articles, while the second includes stochastic chemotherapy scheduling studies. \citet{lameetal16} and \citet{hadidetal22} present a comprehensive review of chemotherapy planning and scheduling studies. We mainly focus on the stochastic chemotherapy scheduling articles, and provide only the most relevant or recent deterministic chemotherapy scheduling articles. Even though there are some similarities between the chemotherapy scheduling and general multi-resource outpatient scheduling studies, we limit the review with articles on chemotherapy processes. We refer the reader to \citet{demiretal21} for the differences between the two streams, and refer to \citet{guptadenton08} and \citet{javidetal17} for a detailed review of studies on other outpatient clinics.  

We start the detailed review with studies considering a deterministic setting for chemotherapy appointments. The fundamental difference between our study and the articles in this category is that we consider uncertainty in infusion durations. \citet{liangturkcan16} and \citet{tabrizietal20} are the most relevant studies of this category to our study. \citet{liangturkcan16} formulate a separate model for the functional and primary care delivery schemes, but do not conduct comparative analysis among them. \citet{tabrizietal20} study a partially flexible setting where each patient is associated with two nurses. Under this setting, the manager does their best to assign patients to their primary nurses or secondary nurses, but also assign some to the floating nurses when necessary. We consider a similar setting to represent the partially flexible delivery scheme, but we assume that all nurses that can be assigned to a patient, except the primary nurses, are given the same priority during the nurse assignment task. In other words, we do not associate a patient with a secondary nurse.  Furthermore, the number of assignments of patients to alternative nurses is controlled through a parameter that is determined by the OCC manager. Therefore, our model can also be adapted to the primary care or functional care delivery schemes.

We next review the most recent deterministic chemotherapy scheduling studies. All models in the following studies are valid for only OCCs that operate according to functional care delivery scheme. \citet{benzaidetal20} study both planning and scheduling phases, and develop three different integer programming (IP) models to assign patients to days, time slots and nurses. \citet{hesarakietal20} formulate a mixed-integer programming model assigning patients to nurses and time slots. \citet{lyonetal22} group treatments of multiple patients in patterns and assigns patterns to time slots in a given day. They schedule also chemotherapy drug preparations in their IP model. The primary objective of the model is to minimize overtime. The secondary objective becomes active when overtime is not used in some patterns to minimize the number of unused slots at the end of the day. \citet{cataldoetal23} consider both planning and scheduling problems. They assign patients to days using a heuristic. They then determine the chairs, nurses, and time slots for the patients using two separate IP models consecutively. The first model is inspired by the model formulated in \citet{turkcanetal12} and minimizes makespan. The second model considers the makespan value as a parameter value and aims to schedule treatments to the earliest possible time slots.  

We next review the studies in the stochastic chemotherapy scheduling category that our article also belongs to. The common difference between our study and the following reviewed ones is related to nursing care flexibility. All models in this category of studies assume functional care delivery scheme. On the other hand, our study analyzes the impact of the nursing care flexibility and proposes a TSMIP model that is applicable to not only functional care delivery, but also primary care and partially flexible care delivery schemes. The relevant studies and the particular differences between each study and our article are discussed next.

\citet{castaingetal16} determine the patient appointment times for a fixed sequence of patients with the objective of minimizing expected waiting time and makespan. All patients are served by the same nurse, therefore only chair assignment decisions are made in the model. Our model is different from \citet{castaingetal16}, since we consider an unfixed patient sequence, multiple nurses, and allow the OCC manager to control nursing care flexibility. \citet{gul21} extends the model in \citet{castaingetal16} by considering multiple nurses and balancing daily workload of nurses. The main difference of our study from \citet{gul21} is that we make patient sequencing decisions and control nursing care flexibility in the TSMIP model.  

\citet{mandelbaumetal20} propose a data-driven approach to schedule infusion appointments by considering chairs as the only resources. However, our study explicitly models nurses and focuses on the impact of nursing care flexibility. \citet{slocumetal20} implement a simple deterministic heuristic to create schedules based on the expected infusion durations and test schedules on a discrete-event simulation model. Since we propose a stochastic optimization model, our study is very different from \citet{slocumetal20}. \citet{alvaradontaimo18} consider both planning and scheduling phases using mean-risk stochastic programming, however the model can schedule one patient at a time. Our model schedules all patients and assigns them to nurses and chairs in a single model. Furthermore, we control the level of nursing care flexibility in the model. \citet{gonzalezmaestroetal22} considers both oncologist and chemotherapy appointments in a single model, but their stochastic model that minimizes patient waiting time is not generic and formulated only for 4 representative scenarios. Furthermore, the stochastic model does not consider nurse assignment decisions. The nurses are assigned in a separate deterministic model. On the other hand, nurse assignment decisions are made by considering uncertain factors in our TSMIP model.

The article that is most relevant to our article is \citet{karakayaetal22}. The authors schedule patient appointments by considering patient acuity levels to minimize expected excess acuity, patient waiting time and nurse overtime. They assume a patient can be assigned to a nurse as long as the skill level of the nurse is larger than the acuity level of the patient. They use a scenario-bundling based decomposition algorithm to find solutions. However, their model is applicable to only functional care delivery scheme. Therefore, they do not study the impact of nursing care flexibility into the performance measures, which is done in our study. Furthermore, they only use random grouping while implementing the scenario-grouping based decomposition algorithm. On the other hand, we propose and test a solution based and two input-based  heuristics to group scenarios while implementing the SGBD algorithm.

\section{Problem Description}
We study a chemotherapy scheduling problem where appointment times are set and patients are assigned to nurses/chairs considering uncertainty in infusion durations under different nursing care delivery schemes. We develop a TSMIP model that is valid for any of the commonly used delivery schemes including primary, functional or partially flexible care. 

The decisions given in the first stage of the TSMIP model are as follows: (i) patient sequencing, (ii) appointment time setting, (iii) patient assignment to nurses, (iv) patient assignment to chairs. We assume that patients are punctual, therefore the patient arrival times are equal to patient appointment times. Patient-to-nurse assignment decisions are made before patients arrive to the OCC on the day of treatment. A patient can be assigned to any nurse that has sufficient skills to treat the patient. However, each patient has a primary nurse. Even though a patient can be assigned to a nurse other than their primary nurse, the total number of patient assignments to alternative nurses may be restricted. The upper limit on the assignment to alternative nurses is controlled by the OCC manager based on the level of flexibility they would like to provide. If the upper limit is set as 0, then the primary care delivery scheme is assumed. If the limit is equal to at least the number of patients, then the functional care delivery is considered. Finally, if the upper limit is in between these two extreme values, then the partially flexible care is in effect. Note that the chairs are identical and therefore patients can be assigned to any of them.  

As is the case in \citet{karakayaetal22}, we assume that the premedication durations are constant for all patients. However, infusion durations are uncertain and their values vary across scenarios. Based on the decisions given in the first stage and the realizations of infusion durations at each scenario, the performance measure values of the study are observed. In other words, no actual decision is given at the second stage of the model. Only patient waiting time and nurse overtime values are recorded for each scenario in this stage. 

Due to the well-known trade-off between patient waiting time and overtime, they are selected as the performance measures of the study \citep{demiretal21}. Furthermore, the perspectives of patients and providers are both taken into account when this trade-off is considered. Patient waiting time is equal to the difference between treatment start time and patient appointment time. A patient may need to wait for a chair or nurse when the treatment of the preceding patient does not finish at the planned time due to uncertainty in infusion durations. Nurse overtime is equal to the difference between the treatment finish time of the last patient that a nurse serves and the shift length of the nurse. Since the number of patient appointments assigned to a day is determined at the chemotherapy planning phase, which is not the focus of this study, nurse overtime may be observed in a given day even if the optimal schedule is implemented. However, nurse overtime may be lowered by carefully setting the patient appointment times considering the uncertainty in infusion durations.    

A nurse can apply only one premedication at a time, since one-to-one relationship between patient and nurse is important during the premedication phase. On the other hand, a nurse may monitor multiple infusions in a given time, since they mainly ensure that patients do not experience any unprecedented event during the infusion activity. As pointed out in \citet{castaingetal16} and  \citet{demiretal21}, patient-to-nurse ratio should not exceed 4 to maintain safe monitoring process. We consider higher patient-to-nurse ratios as unrealistic cases and ignore such parameter settings, therefore we do not explicitly model limit on infusion monitoring processes of nurses. Determining patient-to-nurse ratios is out of the scope of this study, and generally set by considering the practical guidelines. However, the additional constraints can be added as in \citet{demiretal21} in case a formulation that may even apply to unrealistic cases is targeted.    

\begin{figure}
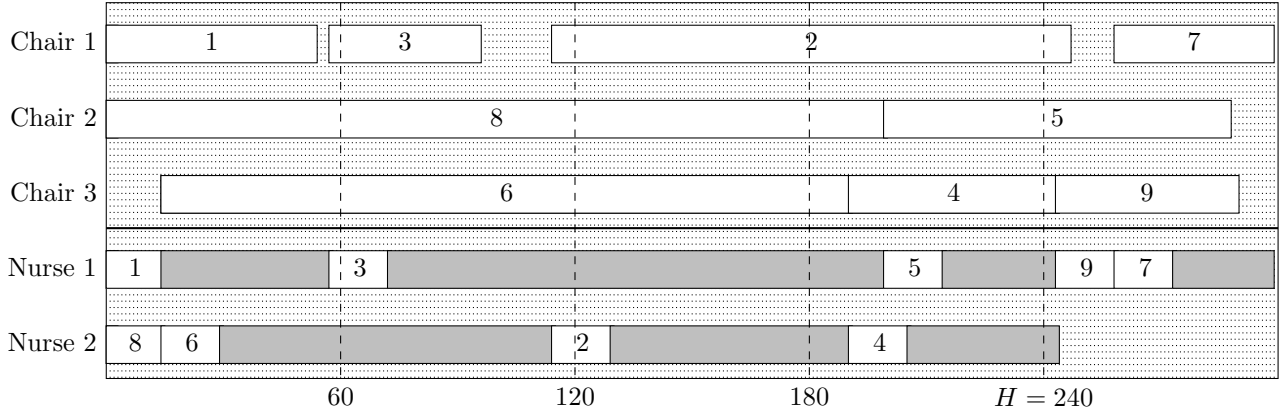

\begin{ganttchart}[vgrid,expand chart=\textwidth, bar height = 0.5]{1}{300}
\ganttbar{Chair 1}{1}{3}
\ganttbar[inline]{1}{1}{54}
\ganttbar[inline]{3}{58}{96}
\ganttbar[inline]{2}{115}{247}
\ganttbar[inline]{7}{259}{299}\\
\ganttbar{Chair 2}{1}{3} 
\ganttbar[inline]{8}{1}{200}
\ganttbar[inline]{5}{200}{288}\\
\ganttbar{Chair 3}{15}{244}
\ganttbar[inline]{6}{15}{191}
\ganttbar[inline]{4}{191}{244} 
\ganttbar[inline]{9}{244}{290}
\ganttnewline[thick]

\ganttbar{Nurse 1}{1}{3}
\ganttbar[inline]{1}{1}{15}
\ganttbar[inline, bar/.append style={fill=gray!50}]{}{15}{58}
\ganttbar[inline]{3}{58}{73}
\ganttbar[inline, bar/.append style={fill=gray!50}]{}{73}{200}
\ganttbar[inline]{5}{200}{215}
\ganttbar[inline, bar/.append style={fill=gray!50}]{}{215}{244}
\ganttbar[inline]{9}{244}{259}
\ganttbar[inline]{7}{259}{274}
\ganttbar[inline, bar/.append style={fill=gray!50}]{}{274}{299}
\\
\ganttbar{Nurse 2}{1}{3}
\ganttbar[inline]{8}{1}{15}
\ganttbar[inline]{6}{15}{30}
\ganttbar[inline, bar/.append style={fill=gray!50}]{}{30}{115}
\ganttbar[inline]{2}{115}{130}
\ganttbar[inline, bar/.append style={fill=gray!50}]{}{130}{191}
\ganttbar[inline]{4}{191}{206}
\ganttbar[inline, bar/.append style={fill=gray!50}]{}{206}{244}

\ganttvrule[vrule/.append style={thin}]{60}{60}
\ganttvrule[vrule/.append style={thin}]{120}{120}
\ganttvrule[vrule/.append style={thin}]{180}{180}
\ganttvrule[vrule/.append style={thin}]{$H=240$}{240}
\end{ganttchart}
\caption{Example of an actual scenario realization for a half-day shift schedule. The numbers represent patients, darker bars for nurses represent the times spent for infusion, while lighter bars for nurses represent premedication durations}
\label{fig:schedule}
\end{figure}

To make the activity flow in the OCC and performance measure definitions clear, we present a scenario realization example for a half-day shift (i.e. 240 minutes) in Figure \ref{fig:schedule}. Suppose that no nurse care flexibility is allowed and therefore assignment to alternative nurses is not possible in this clinic. The figure shows the actual treatment and infusion start times for 9 patients by each of the 2 nurses on each of the 3 chairs. Assume the appointment times for patients 1 through 9 are as follows: $0, 115, 58, 168, 166, 15, 251, 0, 234$. The premedication duration is constant and equal to 15 minutes for all patients, while the infusion durations are 39, 117, 23, 38, 73, 161, 25, 185, 31 minutes, respectively. Furthermore, assume that the primary nurse of the patients with the odd-numbered indices is nurse 1, and that of the patients with the even-numbered indices is nurse 2. Patients 1 and 8 arrive at minute 0 and start their treatment immediately on chairs 1 and 2. Chair 3 stays idle until minute 15, since patient 6 arrives at minute 15 and both nurses are busy with other patients before then. When nurse 2 gets done with the premedication of patient 8, they start the premedication of patient 6. After that premedication ends, nurse 2 only monitors the infusions of patients until the time that the premedication of patient 2 starts. Nurse 1 finishes all premedication and infusion activities at minute 299, and nurse 2 at minute 244.  

In the given example, patient waiting occurs in three cases due to chair unavailability, and in one other case due to nurse care inflexibility leading to nurse unavailability. Patient 5 arrives at minute 166 and waits for 34 minutes until the treatment start, since no chair is available before minute 200. Similarly, patient 4 waits for 23 minutes between times 168 and 191, and patient 9 waits for 10 minutes between times 234 and 244 due to chair unavailability. On the other hand, the waiting of patient 7 occurs due to nurse care inflexibility. When patient 7 arrives at minute 251, they find chair 1 available. However, the patient still waits for 8 minutes, because patient 7's primary nurse (nurse 1) is busy with patient 9 at that time. If the patient was allowed to be assigned to their alternative nurse (nurse 2), the treatment for patient 7 could start immediately upon their arrival as the nurse is available for premedication then. Overtime is observed for both nurses, as nurse 1 works for 59 minutes and nurse 2 for 4 minutes after minute 240.      

Using the notation given in Table \ref{tab:notation} and assuming a finite support for the uncertainty in infusion durations, we formulate the following TSMIP model. 

\begin{table}
\centering
\caption{Notation used in the model}
\label{tab:notation}
\small
\begin{tabular}{cl}
\multicolumn{2}{c}{\textbf{Index sets}} \\
\hline
$P$ & Patients \\ 
$N$ & Nurses \\
$N_i$ & Nurses that can be assigned to patient $i \in P$  \\
$C$ & Chairs	\\ 
$\Omega$ & Scenarios \\

\hline
& \\
\multicolumn{2}{c}{\textbf{Parameters}} \\
\hline
$s$ & Premedication duration of each patient \\ 
$t^{\omega}_{i}$ & Infusion duration of patient $i \in P$ in scenario $\omega \in \Omega$ \\
$H$ & Shift duration of nurses \\
$L$ & Overtime limit of nurses \\
$\lambda$ & Trade-off parameter representing the relationship between waiting time and overtime and taking values in the interval $[0,1]$ \\
$M$ & A large value\\ 
$f_{in}$  & $= \begin{cases} 1,  & \text{if nurse $n \in N$ is the primary nurse of patient $i \in P$} \\0, & \text{otherwise} \end{cases}$\\
$J$ & Upper limit on the number of alternative nurses that can be assigned to patients  \\
\hline
& \\
\multicolumn{2}{c}{\textbf{First-stage decision variables}} \\
\hline
$u_{ij}$ &  $= \begin{cases} 1, & \text{if patient $i \in P$ precedes patient $j \in P$ in the daily appointment list } \\ 0, & \text{otherwise} \end{cases}$ \\
$y_{ic}$  & $= \begin{cases} 1,  & \text{if patient $i \in P$ is assigned to chair $c \in C$} \\0, & \text{otherwise} \end{cases}$\\
$x_{in}$  & $= \begin{cases} 1,  & \text{if patient $i \in P$ is assigned to nurse $n \in N$} \\0, & \text{otherwise} \end{cases}$\\
$a_i$ & Appointment time of patient $i \in P$ \\ 
\hline
& \\
\multicolumn{2}{c}{\textbf{Second-stage decision variables}} \\
\hline

$w^\omega_{i}$ & Waiting time of patient $i \in P$ in scenario $\omega \in \Omega$ \\ 
$o^\omega_{n}$ & Overtime of nurse $n \in N$ in scenario $\omega \in \Omega$\\

\hline
\end{tabular}
\end{table}
\newpage
\begin{flalign}
\mbox{min} \quad & {\mathcal{Q}}(\bm u, \bm x,\bm y, \bm a) & \label{pf1}\\
\mbox{s.t.} \quad & u_{ij}+u_{ji}=1  \hspace{1 cm}  \forall i,j \in P: j>i &\label{pf2}\\
&a_j \geq a_i  -M \big(1-u_{ij}\big) \quad \quad  \forall i, j \in P: j \neq i &\label{pf3}\\
&\sum_{c\in{C}}y_{ic}=1   \hspace{1 cm} \forall i \in P & \label{pf4}\\
&\sum_{n\in{N_i}} x_{in}=1   \hspace{1 cm} \forall i \in P & \label{pf5}\\
&\sum_{i\in I}\sum_{n\in{N_i}} (1-f_{in}) x_{in}  \leq J   \hspace{1 cm}  & \label{pf6}\\
&u_{ij}\in \{0,1\} \hspace{1 cm} \forall i,j \in P: j\neq i & \label{pf7}\\
&y_{ic}\in \{0,1\}  \hspace{1 cm} \forall i\in P,\, \forall c\in C & \label{pf8} \\
&x_{in}\in \{0,1\}  \hspace{1 cm} \forall i\in P,\, \forall n\in N_i & \label{pf9} \\ 
& a_i \geq 0  \hspace{1 cm}  \forall i \in P &\label{pf10} 
\end{flalign}
where $\mathcal{Q}(\bm u,\bm x,\bm y,\bm a)=E_{\xi}\left[Q\left(\bm u,\bm x,\bm y,\bm a,\bm\xi\left(\omega\right)\right)\right]$, and $Q(\bm u,\bm x,\bm y,\bm a,\bm \xi(\omega)))$ represents the following:
\begin{flalign}
\mbox{min}\quad & \lambda\,\sum_{i\in P}w_i^\omega+(1-\lambda)\sum_{n\in N} o_n^\omega& \label{ps1}\\
\mbox{s.t.}\quad 
&a_j + w^\omega_{j} \geq a_i + w^\omega_{i} + s -M \big(3-u_{ij}-x_{in}-x_{jn} \big) \quad \quad \forall i, j\in P,j\neq i,\forall n\in N_i  &\label{ps2}\\ 
&a_j + w^\omega_{j} \geq a_i + w^\omega_{i} + s + t^\omega_{i} -M \big(3-u_{ij}-y_{ic}-y_{jc}\big) \quad \quad \forall i, j \in P,  j \neq i, \forall c \in C & \label{ps3} \\ 
&O^\omega_{n} \geq a_i + w^\omega_{i} + s + t^\omega_{i} -H -M(1-x_{in}) \quad \quad \forall i \in P,\forall n \in N_i  & \label{ps4}\\
&O^\omega_{n} \leq L \quad \quad \forall n \in N &\label{ps5}\\
&w_i^\omega\ge{0} \hspace{1 cm}  \forall{i}\in P &\label{ps6}\\
&O_n^\omega\ge{0}  \hspace{1 cm}  \forall n \in N&\label{ps7}
\end{flalign}

The objective function \eqref{pf1} minimizes the expected total cost of waiting time and nurse overtime. Note that the objective function does not include any first-stage cost. Constraints \ref{pf2} enforce that either patient $i \in P$ precedes $j \in P$  in the daily appointment sequence, or the opposite is true. Constraints \ref{pf3} ensure that the appointment time of patient $j \in P$ is not less than that of patient $i \in P$, if patient $j$ comes after patient $i$ in the appointment sequence. Constraints \ref{pf4} and \ref{pf5} enforce the assignment of patient $i \in P$ to exactly one chair and one nurse, respectively. Constraints \ref{pf6} impose an upper limit, $J$, on the total number of patient assignments to alternative nurses. If $J$ is set to 0, this implies all patients must be assigned to their primary nurses. The case where $J$ is greater than or equal to the total number of patients represent the other extreme. In this case, any patient can be assigned to any alternative nurse. Therefore, the level of flexibility allowed in the system can be controlled through the value of $J$. Constraints \ref{pf7} - \ref{pf9} represent the binary restrictions, and constraints \ref{pf10} represent nonnegativity constraints on the first-stage variables. In the TSMIP model, the appointment times ($a_i$) are represented using continuous variables as is the case in \citet{castaingetal16}. Alternatively, integer variables could be used as in \citet{demiretal21}. We prefer continuous variables, because our solution algorithm works faster and hence larger instances are solved in this case. However, the solution methodology we use is still applicable if the continuous variables are replaced by integer variables.

The second-stage subproblem for each scenario is formulated by \eqref{ps1} - \eqref{ps7}. The objective function \eqref{ps1} minimizes the weighted total cost of patient waiting time and nurse overtime. Constraints \eqref{ps2} ensure that a nurse starts the treatment of patient $j \in P$ after the premedication of patient $i \in P$ ends, in case both patients $i$ and $j$ are assigned to the same nurse, and patient $i$ precedes patient $j$ in the appointment sequence. Constraints \eqref{ps3} ensure that the treatment of patient $j \in P$ starts after the treatment of patient $i \in P$ ends, if patient $i$ precedes patient $j$ in the appointment sequence, and both patients are assigned to the same chair. Overtime amounts for each nurse are calculated through the constraints \eqref{ps4}. Constraints \eqref{ps5} restrict the possible values of nurse overtime from above. Finally, constraints \eqref{ps6} and \eqref{ps7} represent nonnegativity restrictions on the second-stage variables.    

\section{Solution methodology}
Solving the TSMIP model is challenging particularly when the instance sizes are not small. The number of scenarios is naturally one of the important factors affecting the instance size. Even if the model cannot be solved within a reasonable amount of time when formulated for a given scenario set, it may be solved very fast for the case when only a subset of scenarios is considered in the formulation. Furthermore, when the model is solved for multiple disjoint subsets of the original scenario set, a promising solution may be obtained through at least one of them. Based on this observation, we develop and test multiple variants of a decomposition algorithm (SGBD algorithm) that groups scenarios and solves the TSMIP model consecutively based on the created scenario bundles. The variants of the SGBD algorithm differ from each other based on the approach used for grouping scenarios.

The scenario grouping method is frequently used in combination with well-known decomposition algorithms. For example, \citet{crainicetal14}, \citet{escuderoetal13}, \citet{gadeetal16}, \citet{jiangetal21} use it within the progressive hedging algorithm, \citet{escuderoetal13} in a Lagrangian decomposition algorithm, and \citet{oliveiraetal11} in an L-shaped algorithm. Recently, \citet{karakayaetal22} use scenario bundling-based decomposition approach that consists of two steps. They first find the initial solutions for their problem by solving the subproblems created based on disjoint scenario groups. Then, they improve those solutions in the second phase of their approach using an improvement heuristic. Note that they randomly group scenarios in the first phase of their algorithm. Even though random scenario grouping seems to work satisfactorily based on the comparison of the resulting solutions with the optimal solutions \citep{karakayaetal22}, it is worth investigating the following research question: May the solution quality be further improved when efficient algorithms are utilized for grouping? To answer this question, we develop and test multiple SGBD algorithms on the TSMIP model. 

In particular, we test four different methods to partition the original scenario set into multiple scenario groups. The first variant of the SGBD algorithm, called the \textit{progressive SGBD algorithm}, is inspired by the idea used in \citet{hewittetal22} and \citet{keutchayanetal22}. These authors do not group scenarios based on the similarity or dissimilarity of the stochastic input parameter values. They instead group scenarios based on the similarity of single-scenario subproblem solutions. They use sophisticated and time-consuming approaches for this purpose. In particular, \citet{hewittetal22} define graphs on the scenario space and apply graph clustering techniques to group scenarios. \citet{keutchayanetal22} formulate an integer programming model to minimize clustering errors. Both studies test their approaches on the classical stochastic programming models such as the stochastic network design models. To solve the large-size instances of a complex model such as the TSMIP model within a reasonable amount of time, those approaches are not practical and hence more efficient scenario grouping strategies are needed. Therefore, we test a simple solution-based grouping heuristic to create groups progressively. We investigate whether a simple solution-based grouping method may yield better solutions than those of the input-based grouping or random grouping methods when an SGBD algorithm is implemented. As a result, we compare the progressive SGBD with two input-based grouping methods and a random grouping method. In the following subsections, we provide the details of each SGBD algorithm.

\subsection{Progressive SGBD algorithm }

We first present the overall procedures followed throughout the \textit{progressive SGBD algorithm (P-SGBD)}. Initially, the original scenario set having $|\Omega|$ scenarios is divided into $ng = |\Omega|$ groups each having a single element. Next, the TSMIP is solved for each scenario group independently providing $ng$ solutions to be stored. Then, the TSMIP model with the original scenario set, which is called the \textit{master problem}, is solved $ng$ times using the first-stage variable values of each stored solution. This gives $ng$ objective values which are then sorted in a non-decreasing order. Next, every $\alpha$ consecutive groups are combined in a new group starting from the top of the list. This approach yields \(\displaystyle \ceil[\Big] {\frac{ng}{\alpha}}\) groups and finalizes the initial iteration. Note that the last group among the newly created ones may have less than $\alpha$ elements in case \(\displaystyle {\frac{ng}{\alpha}}\) is not an integer. After updating $ng$ as \(\displaystyle \ceil[\Big] {\frac{ng}{\alpha}}\) in the next iteration, the TSMIP is again solved for each group independently. Obviously, the number of models solved decreases, but the time spent for solving each model increases as the iteration count increases. The procedures followed in the first iteration is then repeated to obtain new groups at the end of the second iteration. Then, the same steps are repeated in the following iterations until the iteration counter reaches its limit, $T$. However, the last step (combining scenario groups) is skipped in the last iteration. Finally, the best solution among the solutions obtained for the groups in the last iteration is marked as the SGBD solution. The details of the P-SGBD algorithm is presented in Algorithm \ref{Algo1}.     

\begin{algorithm}
\caption{P-SGBD algorithm}
\label{Algo1}
\begin{algorithmic}[1]
\State \textbf{Step 1:} Initialization
\State {$iter=1, ng=|\Omega|$}
\State Set the iteration limit, $T$
\State Set the number of groups to be combined, {$\alpha$}
\While{$iter \leq T$}
\State {\textbf{Step 2:} Obtain the first-stage decision variable values for scenario group subproblems} 
\For{$g=1,2,\ldots,ng$}
\State Solve scenario group subproblems to find the values of the first-stage variables: $\bm{u}_{g}, \bm{a}_{g}, \bm{x}_{g}, \bm{y}_{g}$
\EndFor
\State{\textbf{Step 3:} Solve the master problem using the first-stage solutions of each group subproblem}
\For{$g=1,2,\ldots,ng$}
\State {Using the values of $\bm{u}_{g}, \bm{a}_{g}, \bm{x}_{g}, \bm{y}_{g}$, solve the master problem and find $z_{g}$}
\EndFor
\State{\textbf{Step 4:} Sort the group subproblems based on their associated master problem objective values in a}
\State{non-decreasing order and place the group indices in the sorted group set, $SG$, accordingly}
\State{$SG=\{g_1,g_2,\ldots,g_{ng}\}$}, where $g_i$ denotes the index of the group with the $i$th best objective value 
\State{$iter \leftarrow iter+1$} 
\If{$iter \leq T$}
\State{\textbf{Step 5:} Starting from the top of the sorted group set, combine every $\alpha$ consecutive groups to obtain} 
\State{$\displaystyle \ceil[\Big] {\frac{ng}{\alpha}}$ new scenario groups}
\State{$ng \leftarrow \displaystyle \ceil[\Big] {\frac{ng}{\alpha}} $} 
\State{$i=1$} 
\While{$i \leq (ng +1 -\alpha)$}
\State{Combine the scenarios in the groups, $g_i,\ldots,g_{i-1+\alpha}$, in a new group}
\State{$i \leftarrow i+\alpha$} 
\EndWhile
\EndIf
\EndWhile
\State{\textbf{Output:} The first-stage solution of the subproblem for group $g_1$ and the associated objective value: $\bm{u}_{g_1}, \bm{a}_{g_1}, \bm{x}_{g_1}, \bm{y}_{g_1}, z_{g_1}$} 
\end{algorithmic}
\end{algorithm}

\subsection{Input-based SGBD algorithms}
We propose and test two different variants of the \textit{input-based SGBD algorithms}, which are inspired by the \textit{k-means clustering} method. We do not incorporate the k-means clustering method into our SGBD algorithm framework, because the size of the groups (clusters) is an output for the method \citep{athurvassilvitskii06}. When some groups have large sizes, the relevant scenario group subproblems may not be solved within a reasonable amount of time as also the results of our preliminary experiments reveal. Therefore, the scenario group size, denoted as $Z$ is predetermined in the input-based SGBD algorithms we test. Both variants of the proposed algorithm choose a reference scenario at each iteration. The first variant groups the scenario with the scenarios that are furthest away from their average (i.e. centroid), while the second variant groups it with the scenarios closest to the centroid. We next summarize the overall procedures related to the each variant. 

\subsubsection{Group with the furthest-SGBD (F-SGBD) algorithm}
In the initialization phase of the F-SGBD algorithm, an arbitrary reference scenario is selected and labeled as the centroid of the first group. In the first step, the distance between each scenario and the centroid is calculated. In particular, the Euclidean distance between the two vectors of infusion durations is measured. Next, the scenario that is the furthest away from the centroid is identified, and added to the group. Then, the centroid of the group is updated by calculating the average infusion durations over all scenarios in the group for each patient. To add new scenarios to the first group, the previous steps are followed at each time until the desired group size is reached. Note that the number of scenarios, $gsize$, may be less than the desired group size, $Z$, in the last group if $\displaystyle \frac{|\Omega|}{Z}$ is not an integer. In the next iteration, starting with another arbitrary reference scenario, the same procedures are followed to create the second group. The iterations then continue until $ng=\displaystyle \ceil[\Big] {\frac{|\Omega|}{Z}}$ groups are generated. 

After creating the scenario groups, the TSMIP model is solved for each group, and $ng$ first-stage solutions are obtained. Then, the master problem (TSMIP with the original scenario set) is solved for each of the $ng$ first-stage solutions. Finally, the lowest master problem objective value and the associated first-stage solution is reported. The complete details of the algorithm can be found in Algorithm \ref{Algo2}.

\begin{algorithm}
\caption{F-SGBD algorithm}
\label{Algo2}
\begin{algorithmic}[1]
\State \textbf{Step 1:} Initialization
\State {Set the desired group size, $Z$}
\State {Set $groupsize=Z, ng=\displaystyle \ceil[\Big] {\frac{|\Omega|}{Z}}, maxdist=0$}
\State Set $grouped (\omega) = false \quad  \forall \omega \in \Omega$

\For{$g = 1,2,\ldots, ng$}
\State Choose an arbitrary $\omega \in \Omega$: $grouped(\omega)=false$
\State{Set grouped $(\omega) = true$. Initialize the set of scenarios included in group $g$: $G_g= \{\omega \}$}
\State Initialize the vector for centroid: $\bm{r}^\omega = \bm{t}^\omega $ 
\If{$(\displaystyle {\frac{|\Omega|}{Z}} \notin \mathbb{Z}^{+}) \, \& \,  (g= ng)$}
\State{ $groupsize = |\Omega| - Z*(ng-1)$}
\EndIf
\For{$z = 1,2,\ldots, groupsize$}
\State \textbf{Step 2:} By calculating the distance between each ungrouped scenario and the centroid of the group, 
\State  identify the scenario that is furthest away from the centroid
\For{$\omega=1,2,\ldots,|\Omega|$}
\If{$grouped(\omega) = false$}
\State Set $dist_\omega$ as the Euclidean distance between scenario $\omega$ and the centroid 

\If {$dist_\omega > maxdist$} 
\State $maxdist  \leftarrow dist_\omega$
\State $ \omega max \leftarrow \omega$
\EndIf
\EndIf
\EndFor
\State \textbf{Step 3:} Add the furthest scenario to group $g$ 
\State $G_g \leftarrow  G_g \cup \omega max$
\State $grouped (\omega max) = true $
\State \textbf{Step 4:} Update the centroid vector
\State $\bm{r}^\omega \leftarrow \displaystyle {\frac {(z*\bm{r}^\omega  + \bm{t}^{\omega max} )} {z+1} }$
\State $maxdist \leftarrow 0$
\EndFor
\EndFor
\State \textbf{Step 5:} Obtain the first-stage decision variable values for scenario  group subproblems
\For{$g=1,2,\ldots,ng$}
\State Solve scenario group subproblems to find the values of the first-stage variables: $\bm{u}_{g}, \bm{a}_{g}, \bm{x}_{g}, \bm{y}_{g}$
\EndFor
\State{\textbf{Step 6:} Solve the master problem using the first-stage solutions of each group subproblem}
\For{$g=1,2,\ldots,ng$}
\State {Using the values of $\bm{u}_{g}, \bm{a}_{g}, \bm{x}_{g}, \bm{y}_{g}$, solve the master problem and find $z_{g}$}
\EndFor
\State{\textbf{Step 7:} Identify the lowest master problem objective value}
\State {Find  $z^* = min_{\{g=1,2,\ldots,ng\}} \{z_g\}$}
\State{\textbf{Output:} Lowest master problem objective value and the associated first-stage solution: $\bm{u}^*, \bm{a}^*, \bm{x}^*, \bm{y}^*, z^*$} 
\end{algorithmic}
\end{algorithm}

\subsubsection{Group with the closest-SGBD (C-SGBD) algorithm}
The only difference between the C-SGBD algorithm and the first variant of the input-based SGBD algorithms is that the scenario that is the closest to the centroid is identified and added to the group rather than the furthest one. The explicit version of the algorithm is provided in Algorithm \ref{Algo3}.  

\begin{algorithm}
\caption{C-SGBD algorithm}
\label{Algo3}
\begin{algorithmic}[1]
\State \textbf{Step 1:} Initialization
\State {Set the desired group size, $Z$}
\State {Set $groupsize=Z, ng=\displaystyle \ceil[\Big] {\frac{|\Omega|}{Z}}, mindist= \infty$}
\State Set $grouped (\omega) = false \quad  \forall \omega \in \Omega$

\For{$g = 1,2,\ldots, ng$}
\State Choose an arbitrary $\omega \in \Omega$: $grouped(\omega)=false$
\State{Set grouped $(\omega) = true$. Initialize the set of scenarios included in group $g$: $G_g= \{\omega \}$}
\State Initialize the vector for centroid: $\bm{r}^\omega = \bm{t}^\omega $ 
\If{$(\displaystyle {\frac{|\Omega|}{Z}} \notin \mathbb{Z}^{+}) \, \& \,  (g= ng)$}
\State{ $groupsize = |\Omega| - Z*(ng-1)$}
\EndIf
\For{$z = 1,2,\ldots, groupsize$}
\State \textbf{Step 2:} By calculating the distance between each ungrouped scenario and the centroid of the group, 
\State  identify the scenario that is the closest to the centroid
\For{$\omega=1,2,\ldots,|\Omega|$}
\If{$grouped(\omega) = false$}
\State Set $dist_\omega$ as the Euclidean distance between scenario $\omega$ and the centroid 

\If {$dist_\omega < mindist$} 
\State $mindist  \leftarrow dist_\omega$
\State $ \omega min \leftarrow \omega$
\EndIf
\EndIf
\EndFor
\State \textbf{Step 3:} Add the closest scenario to group $g$ 
\State $G_g \leftarrow  G_g \cup \omega min$
\State $grouped (\omega min) = true $
\State \textbf{Step 4:} Update the centroid vector
\State $\bm{r}^\omega \leftarrow \displaystyle {\frac {(z*\bm{r}^\omega  + \bm{t}^{\omega min} )} {z+1} }$
\State $mindist \leftarrow \infty$
\EndFor
\EndFor
\State \textbf{Step 5:} Obtain the first-stage decision variable values for scenario  group subproblems
\For{$g=1,2,\ldots,ng$}
\State Solve scenario group subproblems to find the values of the first-stage variables: $\bm{u}_{g}, \bm{a}_{g}, \bm{x}_{g}, \bm{y}_{g}$
\EndFor
\State{\textbf{Step 6:} Solve the master problem using the first-stage solutions of each group subproblem}
\For{$g=1,2,\ldots,ng$}
\State {Using the values of $\bm{u}_{g}, \bm{a}_{g}, \bm{x}_{g}, \bm{y}_{g}$, solve the master problem and find $z_{g}$}
\EndFor
\State{\textbf{Step 7:} Identify the lowest master problem objective value}
\State {Find  $z^* = min_{\{g=1,2,\ldots,ng\}} \{z_g\}$}
\State{\textbf{Output:} Lowest master problem objective value and the associated first-stage solution: $\bm{u}^*, \bm{a}^*, \bm{x}^*, \bm{y}^*, z^*$} 
\end{algorithmic}
\end{algorithm}

\subsection{Random SGBD (R-SGBD) algorithm}
R-SGBD algorithm randomly creates $ng=\displaystyle \ceil[\Big] {\frac{|\Omega|}{Z}}$ scenario groups before solving the subproblems. The remaining parts of the algorithm are the same as those of the input-based algorithms. After the groups are created, the TSMIP model is solved for each group first. Then, the master problem is solved to evaluate $ng$ first-stage solutions. Finally, the master problem objective values are compared to find the minimum one. The algorithm is detailed in Algorithm \ref{Algo4}. 

\begin{algorithm}
\caption{R-SGBD algorithm}
\label{Algo4}
\begin{algorithmic}[1]
\State \textbf{Step 1:} Initialization
\State {Set the desired group size, $Z$}
\State {Set $groupsize=Z, ng=\displaystyle \ceil[\Big] {\frac{|\Omega|}{Z}}$}
\State Set $grouped (\omega) = false \quad  \forall \omega \in \Omega$
\For{$g = 1,2,\ldots, ng$}
\State{Initialize the set of scenarios included in group $g$: $G_g= \emptyset$}
\If{$(\displaystyle {\frac{|\Omega|}{Z}} \notin \mathbb{Z}^{+}) \, \& \,  (g= ng)$}
\State{ $groupsize = |\Omega| - Z*(ng-1)$}
\EndIf
\For{$z = 1,2,\ldots, groupsize$}
\State \textbf{Step 2:} Select an ungrouped scenario randomly and add it to group $g$ 
\State Choose an arbitrary $\omega \in \Omega$: $grouped(\omega)=false$
\State{Set grouped $(\omega) = true$ }
\State $G_g \leftarrow  G_g \cup \omega$
\EndFor
\EndFor
\State \textbf{Step 3:} Obtain the first-stage decision variable values for scenario  group subproblems
\For{$g=1,2,\ldots,ng$}
\State Solve scenario group subproblems to find the values of the first-stage variables: $\bm{u}_{g}, \bm{a}_{g}, \bm{x}_{g}, \bm{y}_{g}$
\EndFor
\State{\textbf{Step 4:} Solve the master problem using the first-stage solutions of each group subproblem}
\For{$g=1,2,\ldots,ng$}
\State {Using the values of $\bm{u}_{g}, \bm{a}_{g}, \bm{x}_{g}, \bm{y}_{g}$, solve the master problem and find $z_{g}$}
\EndFor
\State{\textbf{Step 5:} Identify the lowest master problem objective value}
\State {Find  $z^* = min_{\{g=1,2,\ldots,ng\}} \{z_g\}$}
\State{\textbf{Output:} Lowest master problem objective value and the associated first-stage solution: $\bm{u}^*, \bm{a}^*, \bm{x}^*, \bm{y}^*, z^*$} 
\end{algorithmic}
\end{algorithm}

\section{Experimental study}
We use data from the OCC, based in the Hacettepe University Oncology Hospital in Ankara, Turkey in our experiments \citep{demiretal21}. We first compare four SGBD algorithms with each other to find the best grouping strategy. After determining the variant of the SGBD algorithm to be used in the remaining experiments, we compare the chosen variant with the CPLEX to assess the optimality gap. We next compare the SGBD solutions with the baseline schedules that are used in practice. Then, we investigate the impact of nurse care flexibility into the performance measures. Furthermore, we analyze the impact of the cost coefficients, and number of chairs and nurses. Finally, we calculate and discuss the value of stochastic solutions. 

The SGBD algorithms are implemented on Microsoft Visual C++ 2019 using CPLEX 20.1 Concert Technology. The experiments are performed on a computer with dual Intel Xeon Silver 4114 ten-core processors running at 2.4 GHz and 128 GB of memory. The computational time is limited by 3 hours in the experiments.

\begin{table}
\begin{center}
\centering
\caption {The fraction of patients belonging to each patient type and the interval of infusion durations}
\label{tab:durations}
\begin{tabular}{ |c|c|c| }
 \hline
 \textbf{Patient type} & \textbf{Fraction of patients} & \textbf{Duration interval (mins)} \\   
\hline
   1 & 26.96\% & [16,44] \\ 
   2 & 7.85\% & [29,80] \\ 
   3 & 33.33\% & [74,132] \\
   4 & 31.86\% & [125,217] \\
\hline
\end{tabular}
\end{center}
\end{table}

\begin{table}
\begin{center}
\centering
\caption {Types of patients at each instance of the instance set}
\label{patienttypes}
\begin{tabular}{ |c|c|c|c|c|c|c|c|c|c| }
\cline{2-10}
 \multicolumn{1}{c|}{}& \multicolumn{9}{c|}{\textbf{Patient index}}\\ \hline
 \textbf{Instance \#} & \textbf{1} & \textbf{2} & \textbf{3} & \textbf{4} & \textbf{5} & \textbf{6} & \textbf{7} & \textbf{8} & \textbf{9}  \\   
\hline
1   & 0 & 2  & 0  & 0  & 3  & 2  & 1  & 3  & 2 \\ 
2  & 0 & 3  & 2  &  2 & 1  & 2  & 3  & 3  & 0 \\ 
3   & 2 & 0  &  2 & 2  & 2  & 2  & 3  & 2  & 3\\ 
4   & 3 & 0  & 3  & 0  & 2  & 0  &  3 &  2 & 2\\ 
5   & 3 & 0  & 3  & 2  & 0  & 2  & 2  & 2  & 2\\ 
6  & 3 &  2 &  0 &  2 & 1  & 0  & 3  & 2  & 0\\ 
7   & 0 & 3  & 3  & 2  & 0  & 2  & 2  & 3  & 0\\ 
8   & 1 &  2 &  3 &  3 & 3  & 2  & 3  & 0  & 0\\ 
9   &  2 & 0  & 1  & 2  & 2  & 3  & 0  & 1  & 3\\ 
10   &  0 & 0  & 3  & 2  & 2  & 3  & 0  & 2  & 0\\ 
\hline
\end{tabular}
\label{explicittypes}
\end{center}
\end{table}
\subsection{Problem instances}
We conduct each computational experiment on an instance set with 10 instances. The instances in an instance set differ from each other based on the sampled infusion durations. The number of sampled infusion durations for a patient (i.e., number of scenarios, $|\Omega|$) in each instance is 48 and 96 for small-size instances and moderate-size instances, respectively. The small-size instances are used to compare the SGBD algorithms with each other and the CPLEX, and conduct sensitivity analysis on the number of chairs and nurses, while the moderate-size instances are used for the rest of the experiments. Note that the CPLEX can provide solutions only for small-size instances. The values of the sampled infusion durations depend on the type of patient. Four patient types exist in the data set, which are created by grouping patients based on the treatment durations. The fraction of patients belonging to each patient type, and the interval of infusion durations for each type are provided in Table \ref{tab:durations}, while the additional details can be found in \citet{demiretal21}. The types of patients in an instance are determined in accordance with the fractions given on Table \ref{tab:durations}. The sampled type for each patient at each instance is explicitly shown in Table \ref{explicittypes}. The infusion duration for each patient in each scenario is then sampled from the appropriate time interval using uniform distribution as in \citet{demiretal21}.      

The OCC operates according to two shifts (i.e. morning and afternoon) separated by a lunch break every day. Therefore, we create shift schedules for a time period consisting of 240 minutes (i.e. $H= 240$). We set premedication durations ($s$) as 15 minutes for all patients due to tight 95\% confidence interval for premedication durations in the data set \citep{demiretal21}

We assume that a group of nurses is responsible for only the chairs in a single room since the patients must be monitored continuously. Therefore, we create a schedule for the patients treated in only one of the rooms. In particular, we consider a baseline setting where $|P|=9$ patients can be treated by one of the $|N|=2$ nurses in one of the $|C|=3$ chairs. We induce partial flexibility and allow the assignment of at most 2 patients to their alternative nurses, therefore we set $J=2$ in the baseline instance.  Finally, the cost coefficient of waiting time is set as 0.3 ($\lambda = 0.3$), while that of overtime is 0.7 in the baseline setting, since the data shows that the managers give priority to overtime while creating schedules. Note that we vary the value of $\lambda, J, |N|, |C|$ in the sensitivity analysis experiments. 

Note that when average waiting time and average nurse overtime are explicitly reported for an instance, they represent the averages over scenarios. In particular, the average waiting time represents the average of the sum of the waiting times of all patients over scenarios. Similarly, the average nurse overtime represents the average of the sum of nurse overtimes over scenarios.       

\subsection{Comparison of SGBD algorithms}
In this section, we compare the four variants of the SGBD with each other to determine the best scenario grouping strategy and fix the SGBD algorithm to be used in the later experiments. Before comparing the algorithms, we first perform preliminary analysis to determine important parameter values for them. In particular, we set the desired group size, $Z$, that are used in three of the algorithms: F-SGBD, C-SGBD and R-SGBD. We find that $Z=8$ is an appropriate choice, as the solutions of each SGBD algorithm can be obtained within 3 hours. To have a fair comparison of these algorithms with the P-SGBD, we ensure that the group size reaches 8 in the final iteration of the P-SGBD. Therefore, we perform preliminary analysis to determine the number of groups to be combined, $\alpha$, and the iteration limit, $T$, for the P-SGBD algorithm. Then, we set $\alpha=2$ and $T=4$ based on the trade-off between solution quality and run time. 

\begin{table}
\centering
\caption{Comparison of the four variants of the SGBD with each other in terms of the objective value and run time (in seconds) }
\resizebox{\textwidth}{!}{
\begin{tabular}{|c|c|c|c|c|c|c|c|c|}
\cline{2-9}
 \multicolumn{1}{c|}{}& \multicolumn{2}{c|}{\textbf{P-SGBD}}
 & \multicolumn{2}{c|}{\textbf{F-SGBD}} & \multicolumn{2}{c|}{\textbf{C-SGBD}}
 & \multicolumn{2}{c|}{\textbf{R-SGBD}}\\ \hline
\textbf{Instance \# } &
  \textbf{Time} &
  \textbf{Objective} &
  \textbf{Time} &
  \textbf{Objective} &
  \textbf{Time} &
  \textbf{Objective} &
  \textbf{Time} &
  \textbf{Objective} \\ \hline
1  & 8359.60  & 118.17        & 3018.47 & 113.92 & 2594.60 & 120.06 & 2975.29 & 114.15        \\ 
2  & 7183.20 & 153.49 & 4734.98 & 154.25  & 3387.39& 158.33  & 4384.20 & 154.67 \\ 
3  & 7261.92 & 181.21 & 2445.61 &  180.49 & 4410.90 &   183.51       &  2435.78&  181.87 \\ 
4  & 5253.91  & 132.82 & 2682.28 & 134.67 & 1929.57&    136.10      & 2346.46& 133.63 \\ 
5  & 3090.40 & 158.14 & 1012.70 &  158.95         & 918.92 &   163.95        & 715.22 &         158.77   \\ 
6  & 4668.54 & 108.03 & 1512.03 & 108.38 & 1324.44 & 114.98 & 1455.55 & 109.17 \\ 
7  & 5505.63  & 140.89  & 2778.24  & 138.38  & 2470.02  & 143.65  & 2508.84 & 140.43  \\ 
8  & 3454.00 &  159.86      & 1165.03  & 158.22  & 1788.02  & 158.00  & 1418.15 & 159.71  \\ 
9  & 4779.59  & 127.10         & 1413.07  &   126.56         & 1922.75  & 128.31  & 1479.39 &     130.22      \\ 
10 & 3207.24  &  99.97 & 1217.41  &  100.12 & 995.33  &    104.27      & 1072.42 & 104.25           \\ \hline
\textbf{Average} & 5276.40 & 137.97 & 2197.99 & 137.39 & 2174.19 & 141.12 &  2079.13 & 138.69 \\
\hline
\end{tabular}
}
\label{compSGBD}
\end{table}

In Table \ref{compSGBD}, we report the objective values and run times for 10 instances tested for each SGBD algorithm. The results show that the three variants (P-SGBD, F-SGBD, R-SGBD) yield very similar objective values that are better than that for the C-SGBD. Therefore, we first eliminate C-SGBD algorithm. Among the remaining three variants, F-SGBD and R-SGBD require much shorter computational time than the P-SGBD that is eliminated next. Furthermore, since F-SGBD performs slightly better than the R-SGBD in terms of the objective value, we select F-SGBD for the experiments in the following sections.

The results emphasize that a simple solution-based grouping method cannot outperform input-based grouping or random grouping methods when an SGBD algorithm is implemented. A solution-based grouping method may still be more promising but a more sophisticated and probably time-consuming approach than the P-SGBD is needed. Furthermore, random grouping can compete with input-based grouping in terms of the objective value. If an input-based grouping similar to k-means clustering is to be used, then the one that prioritizes scenarios further away from the centroid while creating groups must be considered.

\subsection{Comparison with the optimal solutions}

\begin{table}
\centering
\caption{Comparison of the F-SGBD with CPLEX in terms of the objective value, solution gap and run time (in seconds)}
\begin{tabular}{|c|c|c|c|c|c|}
\cline{2-6}
 \multicolumn{1}{c|}{}& \multicolumn{3}{c|}{\textbf{Objective}}
 & \multicolumn{2}{c|}{\textbf{Time}} \\ \hline
 \textbf{Instance \# } &
\textbf{CPLEX} & 
  \textbf{F-SGBD} &
  \textbf{Gap (\%)} &
  \textbf{CPLEX} &
  \textbf{F-SGBD}  \\ \hline
1  & 113.15  & 113.92& 0.68& 8214.73 & 3018.47     \\ 
2  & 151.52 & 154.25& 1.80 & 2874.05 & 4734.98 \\ 
3  & 177.67 & 180.49& 1.59 & 6970.19& 2445.61 \\ 
4  & 131.24 & 134.67& 2.61&8560.83 & 2682.28\\ 
5  & 155.10 & 158.95& 2.48& 2140.13& 1012.70 \\ 
6  & 106.53 & 108.38& 1.74&6352.15 & 1512.03\\ 
7  & 137.34  & 138.38& 0.76 & 5936.58& 2778.24\\ 
8  & 155.76  & 158.22& 1.58& 2815.68& 1165.03\\ 
9  & 124.14 & 126.56& 1.95&5338.21 & 1413.07\\ 
10 & 98.11    & 100.12& 2.05& 5014.70&   1217.41    \\ \hline
\textbf{Average} & 135.06 & 137.39 & 1.72 & 5421.73& 2197.98 \\
\hline
\end{tabular}
\label{compcplex}
\end{table}

Table \ref{compcplex} compares the F-SGBD with the CPLEX in terms of the objective value and run time and presents the optimality gap associated with the F-SGBD solutions. The results indicate that the F-SGBD algorithm performs very well and yields near optimal solutions. The average optimality gap is only 1.72 \%, while the range for the gap is between 0.68\% and 2.61 \%. Furthermore, the gap is less than 2 \% in 7 out of 10 instances. The F-SGBD is much more efficient than the CPLEX even when compared on the moderate-size instances, as the CPLEX requires 90 minutes to provide the optimal solution, while the F-SGBD spends 36 minutes on average. 

\begin{table}
\centering
\caption{Comparison of the F-SGBD solutions with the baseline schedules in terms of the objective value, patient waiting time (in minutes) and nurse overtime (in minutes)}
\label{baseline}
\begin{tabular}{|c|c|c|c|c|c|c|}
\cline{2-7}
 \multicolumn{1}{c|}{}& \multicolumn{2}{c|}{\textbf{Objective Value}}
 & \multicolumn{2}{c|}{\textbf{Waiting Time}}  & \multicolumn{2}{c|}{\textbf{Overtime}} \\ \hline
 \textbf{Instance \# } & \textbf{Baseline} & \textbf{F-SGBD} & \textbf{Baseline} & \textbf{F-SGBD} & \textbf{Baseline} & \textbf{F-SGBD} \\ \hline
1 & 414.68  & 117.83 & 716.84 & 35.80 &285.18 & 152.98   \\
2 & 487.64 & 151.66 & 707.48 & 47.70 & 393.43 & 196.21 \\
3  & 507.13 & 181.85 & 703.97 & 27.58 & 422.77  & 247.96 \\
4& 466.77	& 140.04&  708.22 & 43.20 & 363.29 	& 181.54  \\
5&  453.95	&  159.93& 695.76 & 41.66 & 350.32	& 210.62 \\
6&  417.75 &115.48	  & 709.06 & 28.57 & 292.90	& 152.72  \\
7&  473.97 	& 139.56& 727.23 & 36.85 & 365.43	& 183.57 \\
8 & 522.43		& 157.49 & 852.08 & 31.57 & 381.15 &211.46    \\
9&  487.01	 &	127.61 & 778.78 & 34.90 & 361.97 &167.34  \\
10& 459.40		& 101.74 & 806.51 & 37.36 &310.64 &129.32   \\ \hline
\textbf{Average}& 469.07 	& 139.32	& 740.59 & 36.52 & 352.71&  183.37 \\
\hline
\end{tabular}
\end{table}

\subsection{Comparison with the baseline schedules}
To illustrate the benefit of using the F-SGBD in the OCCs, we compare the algorithm with the baseline schedules. The baseline schedules of morning shift are generated by following the rules considered in the OCC whose setting is taken as a reference in this study. In particular, the half of the patients arrive at the first slot at 8:00, while the other half arrive at 10:30, and the regular shift finish time is set as 12:00. The patients with longer expected treatment durations are placed into the first slot, while the remaining are assigned to the second. Ten different baseline schedules are created and compared with the F-SGBD solutions. When a baseline schedule is created based on the rule in practice, the second stage of the TSMIP is solved to find the performance measure values after fixing the first-stage variables. Next, the F-SGBD solution is found for the same instance and compared with the baseline schedule with respect to waiting time, overtime and objective value as shown in Table \ref{baseline}. The results show that the schedules obtained by the F-SGBD algorithm significantly outperform the baseline schedules. The use of fixed-length slots in the baseline schedules particularly leads to very poor performance. An OCC manager may then considerably improve the patient flow in the clinic by ignoring slots and using sophisticated approaches such as the F-SGBD to create daily schedules.        
\begin{table}
\centering
\caption{Sensitivity of average objective value, patient waiting time (in minutes) and nurse overtime (in minutes) to $J$}
\label{flexibility}
\begin{tabular}{|l|c|c|c|}
\hline
\multicolumn{1}{|c|}{} & \textbf{Objective Value} &  \textbf{Waiting Time} & \textbf{Overtime}  \\ \hline
J=0& 145.85 & 50.87 & 172.87  \\
J=1 & 140.99 & 37.41 & 185.37  \\
J=2  & 139.39 & 36.54 & 183.47  \\
J=3& 139.01 &  37.88 & 182.35  \\
J=4 & 139.00 & 37.81  & 182.37  \\\hline
\end{tabular}
\end{table}
\subsection{Impact of nurse care flexibility}
In this section, we test how the nurse care flexibility affects the objective value and performance measures of the study. For this purpose, we vary the value of $J$ by starting from 0, and test each J value by conducting experiments on an instance set consisting of 10 instances. Even though $J$ could take values up to 9, the highest value we test is 4, because the objective value stabilizes from then on. In Table \ref{flexibility}, we observe a diminishing marginal improvement in the average objective value as $J$ increases. When $J$ is increased from 0 to 1, the average objective value improves by 3.33\%.  However, the average objective value decreases only by 1.14\% when $J$ is increased from 1 to 2, and  by 0.28 \% when $J$ is increased from 2 to 3. Providing a little flexibility for nurse care in a non-flexible system (i.e. increasing $J$ from 0 to 1) helps reduce the average waiting time substantially (by 26\%) at the expense of relatively low amount of extra (7\% increase) average overtime. 

The results indicate that the OCCs benefit from increased flexibility in nurse care in terms of efficiency, particularly patient waiting time. On the other hand, allowing a full flexibility in the system is unnecessary. Allowing flexibility even for $1/9$ of the patients may provide substantial advantages. Therefore, the OCCs must prefer a delivery system different from the popular forms including primary care and functional care systems. A small group of patients may be selected from the daily list, for example based on the medical status, and given the chance to be assigned to alternative nurses.    
\begin{table}
\centering
\caption{Sensitivity of average objective value, patient waiting time (in minutes) and nurse overtime (in minutes) to $\lambda$}
\label{lambda}
\begin{tabular}{|l|c|c|c|}
\hline
\multicolumn{1}{|c|}{} & \textbf{Objective Value} &  \textbf{Waiting Time} & \textbf{Overtime}  \\ \hline
$\lambda$=0& 170.58 & 646.96  & 170.58  \\
$\lambda$=0.1 & 165.01 & 85.66 &  173.83 \\
$\lambda$=0.2  &  153.60 &	52.38	& 178.91 \\
$\lambda$=0.3 &  139.32	& 36.52	& 183.37  \\
$\lambda$=0.4 &  123.13	& 23.91	& 189.28 \\
$\lambda$=0.5& 105.58	&16.54&	194.63    \\
$\lambda$=0.6 &   86.98	&11.48	&200.22
  \\
$\lambda$=0.7  & 66.73&	5.96	& 208.54
    \\
$\lambda$=0.8 & 46.00&	3.81&	214.79
  \\
$\lambda$=0.9 & 24.68 &	2.79	& 221.74
   \\
$\lambda$=1 &  0.91	& 0.91	& 1200.00
   \\
\hline
\end{tabular}
\end{table}

\subsection{Sensitivity to cost coefficients}
We evaluate the sensitivity of the cost coefficients to the average objective value, waiting time and nurse overtime by varying the value of $\lambda$ between 0 and 1. Each $\lambda$ value is again tested on an instance set consisting of 10 instances. $\lambda=0$ and  $\lambda=1$ are both unrealistic for the OCC, but they are tested just because they represent the extreme setting. Table \ref{lambda} shows the change in average objective value, waiting time and overtime according to the changes in $\lambda$. As anticipated, the average waiting time decreases and overtime increases as $\lambda$ increases. However, the rate of decrease in waiting time is significantly higher than the rate of increase in overtime. When waiting time is the solely important performance measure ($\lambda=1$), the waiting time drops to almost zero. On the other hand, when overtime is the only relevant measure ($\lambda=0$), the average overtime per nurse would still be 85 minutes.       

\begin{table}
\centering
\caption{Average patient waiting time and nurse overtime values (both in minutes) changing based on different numbers of chairs and nurses. NS means not solved. }
\begin{tabular}{|c|c|c|c|c|c|c|}
\cline{2-7}
 \multicolumn{1}{c|}{}& \multicolumn{3}{c|}{\textbf{Waiting Time}}
 & \multicolumn{3}{c|}{\textbf{Overtime}} \\ \hline
 \textbf{ } &
\textbf{$|N|=1$} & 
  \textbf{$|N|=2$} &
  \textbf{$|N|=3$} &
  \textbf{$|N|=1$} &
  \textbf{$|N|=2$} & 
  \textbf{$|N|=3$} \\ \hline
\textbf{$|C|=3$} &  42.33 & 38.00 & NS & 191.14   & 183.84 &  NS  \\ 
\textbf{$|C|=4$}   & 25.12 & 27.56 & NS & 105.15 & 93.08 & NS \\ 
\textbf{$|C|=5$}   & 17.20 & 19.33 & 18.86 & 59.39 & 48.75 & 46.35 \\ 
\textbf{$|C|=6$}  &  NS   & 10.57  & 9.19 & NS & 18.55  & 16.37\\ \hline

\end{tabular}
\label{nursechair}
\end{table}
\subsection{Impact of the number of chairs and nurses}
We examine the relationship between the performance measures of the study and the number of chairs and nurses. We test $|N|= \{1, 2, 3\}$ as the number of nurses responsible for the patients in the same room. The number of chairs, $|C|$, ranges from 1 to 6 in the experiments. However, we ignore the unrealistic  $(|N|, |C|)$ combinations including  $(1,6)$,   $(3,3)$,  $(3,4)$. 

Table \ref{nursechair} shows the average waiting time and overtime values for different combinations of the number of nurses and chairs, where the average is calculated over ten instances. When the number of chairs increases by keeping the number of nurses constant, both measures improve significantly. When the number of nurses increases while the number of chairs is kept constant, overtime again decreases. However, the impact of the number of nurses on overtime is not as high as the impact of the number of chairs. Furthermore, the waiting time may not even improve as the number of nurses increases due to the trade-off between overtime and waiting time. 

\citet{demiretal21} found that both the number of chairs and nurses are effective to reduce waiting time and overtime at the same time in a fully flexible system. On the other hand, our observations reveal that the number of chairs is more critical to maintain both patient and provider satisfaction simultaneously in a partially flexible system. This is intuitive because the nurses may not be efficiently used in a partially flexible system even if their number increases.

\begin{table}
\centering
\caption{Relative percentage VSS based on the objective values, $z^T$ and $z^{MVP}$, of the F-SGBD and MVP solutions, respectively}
\label{vss}
\begin{tabular}{|l|c|c|c|}
\hline
\multicolumn{1}{|c|}{ \textbf{Instance \# }} & \textbf{F-SGBD} &  \textbf{MVP} & \textbf{VSS (\%)}  \\ \hline
1 &  117.83 &	176.04	& 33.07 \\
2 &  151.66	& 223.03	&  32.00 \\
3 &  181.85 &	261.61	& 30.49 \\
4 &  140.04	& 170.02	& 17.63 \\
5 &  159.93 & 238.70		& 33.00 \\
6 &  115.48	& 145.89	& 20.84 \\
7 &   139.56&	175.00	& 20.25 \\
8 &  157.49	& 179.02	&  12.03 \\
9 &  127.61 &  189.05		&  32.50 \\
10 &  101.74	& 146.31	& 30.46 \\ \hline
\textbf{Average}& 139.32 & 190.47 & 26.85\\
\hline
\end{tabular}
\end{table}

\subsection{Value of stochastic solution}
We estimate the value of stochastic solution (VSS) to show the necessity of representing uncertainty in infusion durations in the TSMIP model. Therefore, we first solve the TSMIP model using F-SGBD algorithm on each of the ten instances and record the objective values $(z^T)$. We then solve the mean value problem (MVP). The MVP is formulated by setting the value of the infusion duration parameter (i.e., $t_i^\omega$) for patient $i \in P$ as the expected infusion duration for the patient in a single-scenario version of the TSMIP model. The first-stage variable values obtained by the MVP solution are then fixed, and the second-stage of the TSMIP model with the original scenario set is solved to obtain the objective value, $z^{MVP}$. The difference, $z^{MVP} - z^T$, yields an estimation of the VSS on each instance.

Table \ref{vss} reports $z^T$ and $z^{MVP}$ values, and the relative percentage VSS  representing the percentage improvement in the objective value achieved by considering uncertainty. The results show that the percentage improvement ranges between 12.03 \% and 33.07 \% across instances, while the average improvement is equal to 26.85 \%. \section{Conclusions}
In this article, we formulate a TSMIP model to create daily schedules for patients, and assign patients to nurses and chairs under uncertainty in infusion durations. The model is unique in the sense that it can be used in an OCC that operate according to any of the commonly used nursing care delivery schemes including primary, functional and flexible care. We minimize the expected patient waiting time and nurse overtime in the objective function. We propose the SGBD algorithm which is a novel decomposition algorithm that utilizes scenario grouping approaches. We compare solution-based and input-based grouping methods with random grouping for the implementation of the SGBD. We select one of the input-based algorithms, F-SGBD, for the rest of the experiments. We compare the SGBD solutions with those of the CPLEX and several baseline schedules. We conduct an extensive sensitivity analysis and evaluate the impact of nursing care flexibility, cost coefficients, and the number of chairs and nurses to the near-optimal schedules. We also estimate the value of the stochastic solution to provide an evidence of the benefit of considering uncertainty in infusion durations. In what follows, we summarize the most important findings of the article. 

Input-based grouping methods may outperform a simple solution-based grouping method while applying the SGBD algorithm. Therefore, to better benefit from the information provided by the solutions when grouping scenarios, a more complex and time consuming but still a practical approach may be necessary. The selected variant, F-SGBD, provides near-optimal solutions and significantly outperforms the baseline schedules used in the OCC. The OCC managers must avoid using fixed-length slots to prevent excessive overtime and waiting time values. 

The OCCs benefit from increased flexibility in nursing care by decreasing patient waiting time substantially. The use of a fully flexible system (functional care delivery) is unnecessary, as allowing a minor level of flexibility still leads to considerable advantages. Therefore, the use of a partially flexible system is the best in terms of system efficiency. To achieve reasonable level of efficiency gains, a small subset of patients may be assigned to nurses other than their primary nurses (i.e., alternative nurses), while the majority are assigned to their primary nurses. Since the continuity of care is an important factor to achieve better medical outcomes, the patients with relatively serious medical conditions may be assigned to their primary nurses. 

In a partially flexible delivery setting, waiting time is more sensitive than nurse overtime to the changes in the cost coefficients of these measures. Furthermore, the number of chairs is more critical than the number of nurses to simultaneously obtain reasonable waiting time and overtime values in a partially flexible system. This finding is different from that for a functional care delivery setting, since \citet{demiretal21} show that both nurses and chairs are effective to reduce overtime and waiting time. Finally, the VSS estimation shows that a significant benefit may be achieved when uncertainty in durations is considered while creating schedules.  

In this article, we only focus on the chemotherapy treatment phase and ignore the patient appointments with oncologists. In a future study, we plan to extend the current model to study the complexities that arise when the OCC managers attempt to maintain coordination between oncologist appointments and chemotherapy treatments that are scheduled to the same day. We also plan to propose a more sophisticated solution-based grouping method rather than a simple rule to group scenarios based on the solutions while implementing the SGBD algorithm to solve a complex model such as the TSMIP. We plan to examine whether the input-based grouping methods would still compete with solution-based methods when more complex approaches are used.

\bibliography{references.bib}

\end{document}